\documentclass[prepint,3p]{elsarticle}
\journal{I don't now}
\usepackage{latexsym,amssymb,amsmath,epsfig}
\def\proof{{\boldmath $Proof.$}\hskip 0.3truecm}

\newtheorem{lm}{Lemma}[section]
\def\proof{{\it Proof.}\nobreak\\}
\def\qed{\hfill$\Box$ \bigskip}
\newtheorem{tm}{Theorem}[section]

\newtheorem{pro}{Proposition}[section]

\newtheorem{open}{Open problem}[section]

\newtheorem{ex}[tm]{Example}

\begin{document}

\title{Mirror bipartite graphs}

\author[UPM]{S.C. L\'{o}pez  \corref{cor1}}\ead{susana@ma4.upc.edu,Tel:0034934134114, Fax:0034934015981}

\author[UPCJG]{F.A. Muntaner-Batle}\ead{famb1es@yahoo.es}

\cortext[cor1]{Corresponding author}

\address[UPM]{ Universitat
Polit\`ecnica de Catalunya. BarcelonaTech, Dept. Matem\`atica Aplicada IV;  Esteve Terrades 5, 08860 Castelldefels, Spain }
\address[UPCJG]{Graph Theory and Applications Research Group.
School of Electrical Engineering and Computer Science.
Faculty of Engineering and Built Environment.
The University of Newcastle.
NSW 2308
Australia }

\begin{abstract}
The concept of mirror bipartite graph appears naturally when studying certain types of products of graphs as for instance the Kronecker product. Motivated by this fact, we study mirror bipartite graphs from the point of view of their degree sequences and of their degree sets. We characterize the sequences of degrees of mirror bipartite graphs.
We also show that from a given set $\mathcal{P}$ of positive integers, we can construct a bipartite graph of
 order $2\max \mathcal{P}$, which is mirror. Furthermore, very little is known for the degree sequences of graphs with loops attached, when the number of loops attached is limited. We show in this paper that mirror bipartite graphs constitute a powerful tool to study the degree sequences of these graphs when the number of loops attached at each vertex is at most $1$.

\vspace{0.2cm}
\end{abstract}

\maketitle

\section{Introduction}

For a bipartite graph $G$, we write $G=(V_1\cup V_2, E)$ to indicate that $V_1$ and $V_2$ are the stable sets of $V(G)$. Define a {\it mirror bipartite graph} of order $2n$ to be a bipartite graph $G=(V_1\cup V_2,E)$ for which there exists a bijective function $\varphi : V_1\rightarrow V_2$ such that, for every pair $u,v\in V_1$, $u\varphi (v)\in E(G)$ if and only if $\varphi (u)v\in E(G)$. We say that $u$ and $\varphi (u)$ are {\it mirror vertices}.
An equivalent statement to this is that $G$ can be drawn in the Cartesian plane in such a way that, the $n$ vertices of $V_1$ are the points $\{(0,0), (0,1),\ldots, (0,n-1)\}$,  the $n$ vertices of $V_2$ are the points $\{(1,0), (1,1),\ldots, (1,n-1)\}$, the edges are straight line segments joining adjacent vertices and the resulting configuration is symmetric with respect to the line $x=1/2$.

Let $G=(V_1\cup V_2, E)$ be a bipartite graph. The {\it bipartite complement} of $G$, denoted by $\bar G^{b}$, is the bipartite graph with stable sets $V_1$ and $V_2$ and edges defined as follows: $xy\in E(\bar G^{b})$ if and only if $xy\notin E(G)$.

\begin{lm}\label{bipertite_complement_mirror}
Let $G$ be a mirror bipartite graph. Then $\bar G^{b}$ is also a mirror bipartite graph.
\end{lm}
\proof Let $G$ be a mirror bipartite graph with stable sets $V_1$ and $V_2$. Then, there exists a bijective function $\varphi: V_1\rightarrow V_2$ such that $u\varphi (v)\in E(G)$ if and only if $\varphi (u)v\in E(G)$, for every $u,v\in V_1$. Thus, $u\varphi (v)\in E(\bar G^b)$ if and only if $\varphi (u)v\in E(\bar G^b)$, for every $u,v\in V_1$. \qed

One of the motivations for introducing mirror bipartite graphs is that they appear when studying certain types of
products as for instance the Kronecker product of graphs. Let $G$ and $H$ be two graphs.
The Kronecker product \cite{WhiRus12} (usually known as direct product) $G\otimes H$ is
the graph with vertex set $V(G)\times V(H)$ and $(a,x)(b,y)\in E(G\otimes H)$ if and only if
$ab\in E(G)$ and $xy\in E(H)$. It is not difficult to check that for an arbitrary graph $H$, when $G=K_2$ and $V(K_2)=\{a,b\}$, the graph
$K_2\otimes H$ is a mirror bipartite graph, with stable sets $\{a\}\times V(H)$ and $\{b\}\times V(H)$ and, for each $x\in V(H)$,  the pair $(a,x), (b,x)$ are mirror vertices. Furthermore, each mirror bipartite graph admits a decomposition of this form. This result appears in the following lemma, which is a particular case of a result found in \cite{LopMun13a}.

Another motivation to study mirror bipartite graphs is due to the fact that with the help of these graphs we will be able to provide an interesting characterization for the degree sequences of $l$-graphs. We call {\it $l$-graphs} these graphs without multiple edges and with at most one loop attached to each vertex.
\begin{lm}\label{lemma_mirror_as_kronecker_product}
Let $G$ be a mirror bipartite graph. Then, there exists a $l$-graph $H$ such that $G\cong H\otimes K_2$.
\end{lm}
\proof Let $G=(V_1\cup V_2, E)$ with $V_1=\{a_1, a_2,\ldots, a_n\}$ and let $\varphi: V_1\rightarrow V_2$ be a bijective function such that $a_i\varphi (a_j)\in E(G)$ if and only if $\varphi (a_i)a_j\in E(G)$, for every $a_i,a_j\in V_1$. We consider a graph $H$ with vertex set $V(H)=V_1$ and edge set $E(H)$ defined by $a_ia_j\in E(H)$ if and only if $a_i\varphi(a_j)\in E(G)$. Let $V(K_2)=\{1,2\}$. Then, the function $f: V_1\times V(K_2)\rightarrow V_1\cup V_2$, defined by $f(a_i, 1)=a_i$ and $f(a_i, 2)=\varphi (a_i)$ is an isomorphism between  $H\otimes K_2$ and $G$. \qed

In this note, we characterize the sequences of degrees of mirror bipartite graphs (Theorem \ref{theo_mirror_iff_bigraphic}). We also show that from a given set $\mathcal{P}$ of positive integers we can construct a bipartite graph of order $2\max \mathcal{P}$, which is mirror (Theorem \ref{theo_sets}) and we completely characterize the degree sequences of $l$-graphs (Theorem \ref{the_new}).

\section{Mirror bigraphic sequences}
For $l$-graphs we define the degree of a vertex to be the number of edges incident with the vertex. That is to say, a loop adds exactly one unit to the degree of the vertex.
Let $P$ be a sequence of nonnegative integers. We say that $P$ is {\it (loop) graphic} if there is a ($l$-)graph $G$ with its vertices having degrees equal to the elements of $P$. In this case, we say that $G$ {\it realizes} $P$. Similarly, let $P$ and $Q$ be two sequences of nonnegative integers. The pair $(P,Q)$ is {\it bigraphic} if there is a bipartite graph $G=(V_1\cup V_2,E)$ with the vertices of $V_1$ having degrees equal to the elements of $P$ and the vertices of $V_2$ having degrees equal to the elements of $Q$. In this case, we say that $G$ {\it realizes the pair} $(P,Q)$. Usually, the elements of $P$ and $Q$ are ordered from biggest to smallest.

The sequence $P$ is {\it mirror bigraphic} if the pair $(P,P)$ is bigraphic and there exists a mirror bipartite graph that realizes $(P,P)$. Our first goal is to prove the following theorem.

\begin{tm}\label{theo_mirror_iff_bigraphic}
The sequence $P$ is mirror bigraphic if and only if the pair $(P,P)$ is bigraphic.
\end{tm}

For the proof of Theorem \ref{theo_mirror_iff_bigraphic}, we will use the following  bipartite version of Havel-Hakimi's theorem \cite{Hakimi,Havel}, see for instance \cite{W}.

\begin{tm}
Suppose $P=(p_1\ge p_2\ge \ldots\ge p_n)$ and $Q=(q_1\ge q_2\ge \ldots\ge q_m)$ are sequences of nonnegative integers. The pair $(P,Q)$ is bigraphic if and only if $(P',Q')$ is bigraphic, where $(P',Q')$ is obtained from $(P,Q)$ by deleting the largest element $p_1$ from $P$ and subtracting $1$ from each of the $p_1$'s largest elements of $Q$.
\end{tm}
Next, we are ready to prove Theorem \ref{theo_mirror_iff_bigraphic}.

\proof By definition, if $P$ is mirror bigraphic then the pair $(P,P)$ is bigraphic. Thus, let us see the converse. We proceed by induction on $n$, where $n$ is the length of the sequence $P$. For $n=1$, the only possible pairs are $((0),(0))$ and $((1),(1))$. These possibilites produce the following two graphs: $G=2K_1$ and $G=K_2$, respectively. Therefore, the statement holds for $n=1$.

Assume now that $(P,P)$ is bigraphic, where $P$ is a sequence of length at most $n$. Then, by induction hypothesis $P$ is also mirror bigraphic. We want to show that if $P'$ is a sequence of length $n+1$, namely $P'=(p'_1\ge p'_2\ge \ldots\ge p'_{n+1})$ and $(P',P')$ is bigraphic then, $P'$ is mirror bigraphic. Since $(P',P')$ is bigraphic, it follows that we can apply Havel-Hakimi's theorem twice obtaining a bigraphic pair of identical sequences $(P'',P'')$, where $P''$ is obtained from $P'$ by eliminating $p_1'$ and subtracting $1$ to each element of $p_2', p_3',\ldots, p'_{p_1}$. The rest of the elements remain the same. By induction hypothesis, $P''$ is mirror bigraphic. Let $G''$ be a mirror bipartite graph that realizes $(P'',P'')$. Let $\{a_2,b_2\}, \{a_3,b_3\}, \ldots, \{a_{n+1},b_{n+1}\}$ be the pairs of mirror vertices in $G''$, where $deg_{G''}(a_i)=deg_{G''}(b_i)=p_i''$, for every $i=2,3,\ldots, n+1$. Then, adding a new vertex to each stable set of $V(G'')$, namely $a_1$ and $b_1$ and joining the vertices $a_1$ and $b_1$ by an edge, and all vertices of the form $b_i$ to $a_1$ and the vertices of the form $a_i$ to $b_1$, for every $i=2,3,\ldots, p_1'$, we are done.\qed

\begin{tm}\label{the_new}
Let $P$ be a sequence of nonnegative integers. Then the following statements are equivalent.
\begin{itemize}
  \item[(i)] $P$ is loop graphic.
  \item[(ii)] $(P,P)$ is bigraphic.
  \item[(iii)] $P$ is mirror bigraphic.
\end{itemize}
\end{tm}
\proof By Theorem \ref{theo_mirror_iff_bigraphic}, conditions (ii) and (iii) are equivalent. We will prove that (i) implies (iii) and viceversa. Suppose there exists a $l$-graph that realizes $P$. Then, $L\otimes K_2$ is a mirror bipartite graph that realizes $(P,P)$. Hence, by definition, $P$ is mirror bigraphic. Suppose now that $P$ is mirror bigraphic. Then, by definition there exists a mirror bipartite graph $G$ that realizes $(P,P)$. Thus, by Lemma \ref{lemma_mirror_as_kronecker_product}, there exists a $l$-graph $H$ such that $G\cong H\otimes K_2$. Hence, $H$ realizes $P$.
\qed

Next, let $P=(p_1\ge p_2\ge \ldots\ge p_n)$ be a sequence of nonnegative integers for which $(P,P)$ is bigraphic. Let  Bipp$(P,P)$ be the bipartite graphs (modulo isomorphism) that realize $(P,P)$, and Mirr$(P,P)$ be the mirror bipartite graphs (modulo isomorphism) that realize $(P,P)$.

It is clear that, for each $n\in \mathbb{N}$, the sequence of length $n$, $P_n=(n-1\ge n-1\ge \ldots\ge n-1)$ is a sequence for which Bipp$(P_n,P_n)$=Mirr$(P_n,P_n)$. The next lemma introduces another sequence $P_n$ for which all bipartite graphs that realizes the pair $(P_n,P_n)$ are mirror bipartite graphs.

\begin{lm}\label{lemma_mirror_sequence}
Let $P_n=(n\ge n-1\ge \ldots\ge 1)$, for each $n\in \mathbb{N}$. Then, Bipp$(P_n,P_n)$=Mirr$(P_n,P_n)$.
\end{lm}

\proof It is easy to show that the pair $(P_n,P_n)$ is bigraphic for every $n\in \mathbb{N}$. Next, we will show that,  for every $n\in \mathbb{N}$, there exists a unique bipartite graph (modulo isomorphisms) that realizes $(P_n,P_n)$. Suppose that $G=(V_1\cup V_2, E)$ is a bipartite graph that realizes $(P_n,P_n)$, where $V_1=\{a_i:\ i=1,2,\ldots, n\}$ and deg$_G(a_i)=n-i+1$. Let $V_2=\{b_i:\ i=1,2,\ldots, n\}$. Then, $a_1$ is adjacent to every vertex of $V_2$. Since $1\in P_n$, one of the vertices in $V_2$ should have degree $1$. Without loss of restriction, assume that deg$_G(b_n)=1$. Then,  $a_2$ should be adjacent to every vertex of $V_2\setminus \{b_n\}$. Since $2\in P_n$, one of the vertices in $V_2$ should have degree $2$. Without loss of restriction, assume that deg$_G(b_{n-1})=2$. Then, the vertex $a_3$ should be adjacent to every vertex of $V_2\setminus \{b_{n-1},b_n\}$. We proceed in this way until we complete the adjacencies of all vertices in $V_1$.\qed

\begin{figure}[ht]
\begin{center}
  \includegraphics[width=104pt]{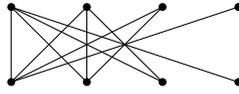}\\
  \caption{The bipartite graph that realizes $(P_4,P_4)$}\label{Fig_1}
\end{center}
\end{figure}

The bipartite graph (modulo isomorphisms) that realizes $(P_4,P_4)$ is shown in Figure \ref{Fig_1}.

It is trivial that $0$-regular graphs of even order and graphs of the form $nK_2$ are mirror bipartite graphs. The following lemma is also easy to prove.

\begin{lm}\label{mirror_2_regular}
Every $2$-regular bipartite graph is a mirror bipartite graph.
\end{lm}
\proof Since every $2$-regular graph is the disjoint union of cycles, it suffices to prove that every cycle $C$ of even order $2n$ is mirror. Let $V(C)=\{v_i\}_{i=0}^{2n-1}$, $E(C)=\{v_iv_{i+1}\}_{i=0}^{2n-2}\cup \{v_0v_{2n-1}\}$, $V_1=\{v_{2i}\}_{i=0}^{n-1}$ and $V_2=\{v_{2i+1}\}_{i=0}^{n-1}$ Then, clearly $V(C)=V_1\cup V_2$ and $\varphi (v_{2i})=v_{2n-1-2i}$, for $i=0,1,\ldots, n-1$, is a function from $V_1$ to $V_2$ such that $v_{2i}\varphi(v_{2j})\in E(G)$ if and only if $\varphi(v_{2i})v_{2j}\in E(G)$, for all $i,j\in \{0,1,\ldots, n-1\}$.\qed

At this point, we are ready to state and prove the following proposition.
\begin{pro}
Let $G=(V_1\cup V_2,E)$ be a bipartite regular graph with $|V_1|=|V_2|=n$. If $G$ is not a mirror bipartite graph then $n\ge 6$.
\end{pro}
\proof It is clear that $n\ge 3$. By previous comments and Lemma \ref{mirror_2_regular}, if $n=3$ then $G$ cannot be neither $0$, $1$, nor $2$-regular. Thus, $G$ is $3$-regular. But in this case, $G\cong K_{3,3}$ and $K_{3,3}$ is mirror. If $n=4$ or $n=5$, then $G$ is $r$-regular, with $0\le r\le 5$. Again, it is not possible for $G$ to be either $0$, $1$ or $2$-regular. However, if $G$ is $r$-regular, with $3\le r\le 5$, then the bipartite complement of $G$ is $r$-regular, with $0\le r\le 2$, and hence, mirror. Therefore, by Lemma \ref{bipertite_complement_mirror}, $G$ is also mirror, a contradiction. This implies that $n\ge 6$.\qed

\begin{ex}
Figure 2 shows a $3$-regular non-mirror bipartite graph of order $12$. This fact is clear since the vertices $u$ and $v$ that appear in one of the stable sets are twin vertices (they share the same set of neighbors), whereas in the other stable set there are not twin vertices.
\begin{figure}[ht]
\begin{center}
  \includegraphics[width=161pt]{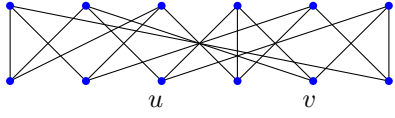}\\
  \caption{A bipartite graph which is not mirror.}\label{Fig_2}
\end{center}
\end{figure}
\end{ex}

\begin{open}
Characterize the sequences $P_n$ of length $n$ for which Bipp$(P_n,P_n)$=Mirr$(P_n,P_n)$.
\end{open}

\section{Bigraphic sets}
Let $\mathcal{P}=\{p_1>p_2>\ldots >p_k\}$ be a set of positive integers. We say that $\mathcal{P}$ is a graphic set if there exists a sequence of the form $P=(p_1\ge \ldots\ge p_1\ge p_2\ge \ldots\ge p_2\ge \ldots \ge p_k\ge\ldots\ge p_k)$, which is graphic. It is an easy observation that every set of positive integers is a graphic set. However, Kapoor et al. introduced in \cite{Kapoor} the following result. A short proof of it can be found in \cite{TriVij07}.

\begin{tm}\cite{Kapoor} \label{Kapoor_theorem}
Let $\mathcal{P}=\{p_1>p_2>\ldots >p_k\}$ be a set of positive integers. Then there exists a graph of order $p_1+1$ with degree set $\mathcal{P}$.
\end{tm}


Let $\mathcal{P}=\{p_1>p_2>\ldots >p_{k_1}\}$ and $\mathcal{Q}=\{q_1>q_2>\ldots >q_{k_2}\}$ be two sets of integers. We say that the pair $(\mathcal{P},\mathcal{Q})$ is bigraphic if there exists a bipartite graph $G$ with stable sets $V_1$ and $V_2$ such that the vertices of $V_1$ have degrees $(p_1\ge\ldots\ge p_1\ge p_2\ge\ldots\ge p_2\ge \ldots\ge p_{k_1}\ge\ldots\ge p_{k_1})$ and the vertices of $V_2$  have degrees $(q_1\ge\ldots\ge q_1\ge q_2\ge\ldots\ge q_2\ge \ldots\ge q_{k_2}\ge\ldots\ge q_{k_2})$. We say that $G$ {\it realizes the pair} $(\mathcal{P},\mathcal{Q})$. Next, we have the following theorem.

\begin{tm}\label{theo_sets}
Let $\mathcal{P}=\{p_1>p_2>\ldots >p_{k}\}$ be a set of integers. Then, $(\mathcal{P},\mathcal{P})$ is bigraphic and there exists a bipartite graph that realizes the pair, with each stable set of size $p_1$. Furthermore, such a graph can be chosen to be a mirror bipartite graph.
\end{tm}

\proof In order to prove the theorem, we will consider two cases.

\underline{Case 1.} Assume that all elements of $\mathcal{P}=\{p_1>p_2>\ldots >p_{k}\}$ form a set of consecutive integers.
 That is to say,  $\mathcal{P}=\{l>l-1>\ldots >l-k+1\}$. In this case,
 we consider the sequence $P'$ obtained from $\mathcal{P}$ in such a way that all elements of $P'$
 are the elements of $\mathcal{P}$ minus $p_k-1$, that is $P'=(k>k-1>\ldots >1)$.
 By Lemma \ref{lemma_mirror_sequence} we obtain that the sequence $P'$ is mirror bigraphic,
 and hence there is a mirror bipartite graph that realizes $(P',P')$, with each stable set of size $k$. Let such a graph be $H$ and let $A$ and $B$ be the stable sets of $V(H)$. Add a new vertex to $A$ and a new vertex to $B$, namely $a_1$ and $b_1$, and join $a_1$ with all vertices of $B$, including $b_1$. Similarly, join $b_1$ with all vertices of $A$. In this way, we obtain a new bipartite graph $H_1$ with vertex set, $A_1=A\cup \{a_1\}$ and $B_1=B\cup \{b_1\}$ such that the vertices of $A_1$ and of $B_1$ have the following degree sequence:
$k+1, k+1,k,k-1,\ldots, 2.$
From $H_1$ we can obtain a new graph $H_2$ with stable sets $A_2$ and $B_2$, in a similar way,
where the vertices of $A_2$ and of $B_2$, have degree sequence: $k+2,k+2, k+2,k+1,k,\ldots, 3.$
We proceed in this way until we reach the graph with the required degree set. This concludes case 1.

\underline{Case 2.} Assume that the elements of $\mathcal{P}=\{p_1>p_2>\ldots >p_{k}\}$ do not form a set of consecutive integers.  Then, there exists $j\in \{1,2,\ldots, k\}$ such that $p_j-p_{j+1}>1$. Let $i$ be the smallest such $j$. Consider the new set $\mathcal{P}'=\{p_1-1>p_2-1>\ldots >p_i-1>p_{i+1}>p_{i+2}>\ldots>p_{k}\}$.  By Theorem \ref{Kapoor_theorem}, there is a graph that has order $p_1$ and degree set $\mathcal{P}'$. Call such a graph $G'$. Assume that the vertices of $G'$ are $a_1,a_2,\ldots, a_{p_1}$ and consider the graph $G''=G'\otimes K_2$. By construction, this new graph is a mirror bipartite graph of order $2p_1$ and with degree sequences in each stable set equal to the degree sequence of $G'$. Moreover, we can add to $E(G'')$ a set of edges of the form $\{((a_j,1),(a_j,2)): \ j=1,2,\ldots, i\}$, where $V(K_2)=\{1,2\}$, in order to obtain a graph $G$ that realizes the pair $(\mathcal{P},\mathcal{P})$. \qed

\section{Conclusions}
We have introduced the concept of mirror bipartite graphs that naturally appears when studying certain types of products of graphs as for instance the Kronecker product. We have studied mirror bipartite graphs from the point
of view of their degree sequences and of their degree sets. The main contributions of this note are Theorem \ref{theo_mirror_iff_bigraphic}, Theorem \ref{the_new} and Theorem \ref{theo_sets}.
Theorem \ref{theo_mirror_iff_bigraphic} establishes that, for a given sequence of positive integers $P$, the pair $(P,P)$ is bigraphic if and only if there exists a mirror bipartite graph that realizes $(P,P)$. Theorem \ref{the_new} presents a characterization of loop graphic sequences, in terms of bigraphic sequences and of mirror bigraphic sequences. In fact, mirror bigraphic sequences can be thought as the link between bigraphic sequences and loop graphic sequences. In Theorem \ref{theo_sets}, we prove that, for a given set of positive integers $\mathcal{P}$, the pair $(\mathcal{P},\mathcal{P})$ can be realized by a  bipartite graph of minimum order, that is, $2\max \mathcal{P}$, which in fact can be chosen to be a mirror bipartite graph.

\noindent {\bf Acknowledgements} The research conducted in this document by the first author has been supported by the Spanish Research Council under project
MTM2011-28800-C02-01 and  by the Catalan Research Council
under grant 2009SGR1387.

\end{document}